\documentclass[a4paper,11pt,twoside]{article}
\usepackage{pst-fill,pst-grad}
\usepackage{textcomp} 
\usepackage[english]{babel}  
\usepackage[utf8]{inputenc}    
\usepackage[titletoc]{appendix}
\usepackage{titlesec}  
\usepackage{graphicx} 
\usepackage{amsmath}  
\usepackage{float} 
\usepackage{fancyhdr}  
\usepackage[matrix,arrow,curve]{xy} 
\usepackage{pstricks} 
\usepackage{amsmath,amsfonts,verbatim,afterpage,theorem,euscript,mathrsfs,amssymb}
\usepackage{amsfonts}
\usepackage{amssymb} 
\usepackage{array} 
\usepackage{dsfont}
\usepackage[colorlinks=true,linkcolor=blue,citecolor=red]{hyperref}
\usepackage{authblk}
\usepackage{color}
 
\newtheorem{Proposition}{Proposition}[section]

\newtheorem{Theoreme}{Theorem}[section]

\def \vu{\vec{u}}
\def \vb{\vec{b}}
\def \ve{\vec{e}}

\def \R{\mathbb{R}}

\def \Rt{\mathbb{R}^3}

\def \finpv{\hfill $\blacksquare$} 
\def \pv{{\bf{Proof.}}~} 

\def \ds{\displaystyle}

\title{\bf A general Liouville-type theorem for the 3D steady-state Magnetic-Bénard system}

\author[1]{ Oscar Jarr\'in\footnote{corresponding author: oscar.jarrin@udla.edu.ec}}

\affil[1]{\scriptsize Escuela de Ciencias Físicas y Matemáticas, Universidad de Las Américas, Vía a Nayón, C.P.170124, Quito, Ecuador.} 
\date{\today}
\begin{document}
	\maketitle	 
\begin{abstract}
We establish a Liouville-type theorem for the elliptic and incompressible Magnetic-Bénard system defined over the entire three-dimensional space. Specifically, we demonstrate the uniqueness of trivial solutions under the condition that they belong to certain local Morrey spaces. Our results generalize in two key directions: firstly, the Magnetic-Bénard system encompasses other significant coupled systems for which the Liouville problem has not been previously studied, including the Boussinesq system, the MHD-Boussinesq system, and the Bénard system. Secondly, by employing local Morrey spaces, our theorem applies to Lebesgue spaces, Lorentz spaces, Morrey spaces, and certain weighted-Lebesgue spaces. \\[3mm]
\textbf{Keywords} 3D Magnetic-Bénard system; Liouville-type theorem; Local Morrey spaces; Caccioppoli-type estimates. \\[2mm]
\textbf{AMS Classification:} Primary 35B53,	Secondary 35B45. 
\end{abstract}

\section{Introduction} 
This short note is devoted to studying uniqueness issues of the steady-state (time-independent) and incompressible 3D Magnetic-Bérnad system defined on the whole space $\Rt$. For a velocity $\vu=\vu(x)$, a scalar pressure $P(x)$, a magnetic vector field $\vb(x)$, a scalar temperature $\theta(x)$, and $\ve_3$ the third vector of the canonical base in $\Rt$,  this system is written  as follows:
\begin{equation}\label{System}
    \begin{cases}\vspace{2mm}
    -\nu \Delta \vu + (\vu \cdot \vec{\nabla})\vu - (\vb \cdot \vec{\nabla}) \vb + \vec{\nabla} \left(P +\frac{1}{2}| \vb|^2 \right) - \zeta \theta \ve_3=0,  \\ \vspace{2mm}
    -\mu \Delta \vb + (\vu \cdot \vec{\nabla})\vb - (\vb \cdot \vec{\nabla})\vu =0,  \\ \vspace{2mm}
    -\kappa\Delta \theta + \vu \cdot \vec{\nabla}\theta - \eta \vu \cdot \ve_3=0, \\
    \text{div}(\vu)=\text{div}(\vb)=0.    
    \end{cases}
\end{equation}

Physically, this system models  the heat convection phenomenon influenced by velocity, magnetic field and temperature \cite{Mulone,Nakamura}, where $\nu>0$ and $\mu>0$ are the viscosity and diffusivity, meanwhile $\kappa>0$ is the thermal diffusion. In this first equation of (\ref{System}) the term $\theta \ve_3$ represents  the acting buoyancy force on the fluid motion with a constant $\zeta >0$. Moreover, in the third equation of (\ref{System}), the term  $\eta \vu \cdot \ve_3$ models the Rayleigh–Bénard convection in a heated inviscid fluid with  a constant $\eta\geq 0$. 

\medskip

Mathematically, the coupled system (\ref{System}) is the time-independent counterpart of the parabolic Magnetic-Bénard system:
\begin{equation}\label{System-Parabolic}
    \begin{cases}\vspace{2mm}
    \partial_t \vu -\nu \Delta \vu + (\vu \cdot \vec{\nabla})\vu - (\vb \cdot \vec{\nabla}) \vb + \vec{\nabla} \left(P +\frac{1}{2}| \vb|^2 \right) - \zeta\theta \ve_3=0, \quad \text{div}(\vu)=\text{div}(\vb)=0, \\
    \partial_t \vb -\mu \Delta \vb + (\vu \cdot \vec{\nabla})\vb - (\vb \cdot \vec{\nabla})\vu =0, \quad    \quad 
    \partial_t \theta -\kappa\Delta \theta + \vu \cdot \vec{\nabla}\theta - \eta \vu \cdot \ve_3=0,
    \end{cases}
\end{equation}
which  has recently garnered  interest in the research community in the analysis of Partial Differential Equations modeling fluid dynamics. Initial studies focused on investigating smooth solutions to  (\ref{System-Parabolic}). In \cite{Nakamura2} the authors proven  smoothness and  time-analyticity of  strong solutions. In \cite{Cheng} it is considered the Cauchy problem in the 2D case  with mixed partial viscosity, and  in  \cite{Zhu}  this work was extended to the  zero thermal conductivity setting. The 3D case was studied in \cite{Ma}, where the authors proven the existence of smooth solutions with mixed partial dissipation, magnetic diffusion and thermal diffusivity, which arise  from small $H^1$-initial data. On the other hand, uniqueness issues for the problem  (\ref{System-Parabolic})  were studied in \cite{Miao}. Precisely, for a constant $A>0$ and for a time $T>0$,  it is proven that this problem  has at most one smooth  solution on $[0,T]\times \Rt$, provided that $|\vec{\nabla}\otimes \vu(t,x)| + |\vec{\nabla}\otimes \vb(t,x)| + |\vec{\nabla}\theta(t,x)| \leq A$ for all $x\in \Rt$, and that  $| P(t,x)| \lesssim |x|^{-\frac{1}{2}}$ when $|x|\to +\infty$. 

\medskip

Coming back to the steady-state system (\ref{System}),  to our knowledge, this system   has not been studied  before, and following the ideas of \cite{Miao}, we aim to  prove some \emph{sufficient} conditions to ensure the uniqueness of classical solutions.  
With a minor loss of generality,  we will set all the physical constants $\nu=\mu=\kappa=\zeta=1$ since they do not play any substantial role in our study. In contrast, in Rayleigh–Bénard convection term  $\eta \vu \cdot \ve_3$,  we shall keep the constant $\eta\geq 0$ in order to consider both active effects (when $\eta>0$) and inactive effects (when $\eta=0$) of this term. We shall see that these cases  become relevant when studying uniqueness issues of solutions to the system (\ref{System}). 

\medskip

Recall that a classical solution to the system (\ref{System}) is the set of functions  $\vu, \vb, \theta \in \mathcal{C}^{2}(\Rt)$ and $P \in \mathcal{C}^1(\Rt)$, which verifies (\ref{System}) pointwise in $\Rt$. However, our first key remark is that   in this general framework, the system (\ref{System}) admits several solutions. For instance, the trivial solution $(\vu,\vb,\theta,P)=(0,0,0,0)$ and also a specific non-trivial solution  defined for any  $x=(x_1,x_2,x_2)\in \Rt$ as follows:
\begin{equation}\label{Non-Trivial-Solution}
\begin{cases}\vspace{2mm}
    \vu(x)= \vb(x)=\left(ax_1,bx_2,-2cx_3\right), \  \theta(x)=4c^2x_3,  \ P(x)=-\frac{a^2}{2}x^{2}_{1}-\frac{b^2}{2}x^{2}_{2}, \\
a+b=2c,  \ \  c=\pm \frac{1}{2} \ \ \mbox{and} \  \ \eta=1.    
    \end{cases}
\end{equation}
Using basic rules of vector calculus, one can verify that these functions satisfy each equation in (\ref{System}). Refer to Appendix \ref{Appendix} for the details. In this way, we seek additional \emph{a priori} conditions on classical solutions to ensure the uniqueness of the trivial solution. This problem is known as the Liouville-type problem for the elliptic system (\ref{System}).

\medskip

The Liouville-type problem has been extensively studied for two particular cases of the general coupled system (\ref{System}). First, setting $\vb=\theta=0$ and $\eta=0$ gives the well-known steady-state Navier-Stokes equations:
\begin{equation}\label{NS}
-\Delta \vu + (\vu \cdot \vec{\nabla}) \vu + \vec{\nabla}P=0, \quad \text{div}(\vu)=0.
\end{equation}
The celebrated result \cite[Theorem X.9.5]{Galdi} shows that if $\vu \in L^{\frac{9}{2}}(\mathbb{R}^3)$, then $\vu=0$ and $P=0$. This result was extended in \cite{Kozono} to the setting of Lorentz spaces, and further improved in \cite{Jarrin2}. For more interesting results, refer to \cite{Chae1,Chae2,Chamorro,Seregin} and the references therein.

\medskip

On the other hand, setting $\theta=0$ and $\eta=0$, the coupled system reduces to the Magneto-hydrodynamic (MHD) equations:
\begin{equation}\label{MHD}
\begin{cases}\vspace{2mm}
- \Delta \vu + (\vu \cdot \vec{\nabla})\vu - (\vb \cdot \vec{\nabla}) \vb + \vec{\nabla} \left(P +\frac{1}{2}| \vb|^2 \right)=0, \quad \text{div}(\vu)=0, \\
-\Delta \vb + (\vu \cdot \vec{\nabla})\vb - (\vb \cdot \vec{\nabla})\vu =0, \quad \text{div}(\vb)=0,
\end{cases}
\end{equation}
for which there is also a significant amount of literature on the Liouville-type problem. See, for instance, \cite{Chen,Cho,Li,Li2,Schulz}.
\medskip

Returning to the coupled system (\ref{System}), in the non-trivial solution given in (\ref{Non-Trivial-Solution}), we observe that 
\[ \left|\big(\vu(x),\vb(x),\theta(x),P(x)\big)\right| \to +\infty, \quad \mbox{as} \quad  |x| \to +\infty.\]
This observation, combined with the previous studies on equations (\ref{NS}) and (\ref{MHD}), strongly suggests to consider functional spaces that characterize \emph{decaying properties} of functions. These spaces are crucial for addressing the Liouville-type problem for (\ref{System}). 

\medskip

There is a wide range of these spaces in the literature, including the classical Lebesgue spaces and the more sophisticated Lorentz and Morrey spaces. Our goal is to work within a sufficiently general framework that encompasses these spaces and others. Inspired by our previous work \cite{Jarrin1}, our framework will be defined by the following \emph{local Morrey spaces}: for $\gamma \geq 0$ and $1 < p < +\infty$, we define $M^{p}_{\gamma}(\mathbb{R}^3)$ as the space of locally $L^p$-integrable functions such that
\begin{equation}\label{Local-Morrey}
\| f \|_{M^{p}_\gamma} := \sup_{R \geq 1} R^{-\frac{\gamma}{p}} \left( \int_{|x| < R} |f(x)|^p  dx \right)^{\frac{1}{p}} < +\infty.
\end{equation}
Here, $\gamma \geq 0$ characterizes the asymptotic behavior of the local quantity $\ds{\left( \int_{|x| < R} |f(x)|^p  dx \right)^{\frac{1}{p}}}$ as $R \to +\infty$. Additionally, we define the space $M^{p}_{\gamma,0}(\mathbb{R}^3)$ as the subset of $M^{p}_{\gamma}(\mathbb{R}^3)$ consisting of functions $f$ that satisfy the additional vanishing property at infinity:
\begin{equation}\label{Local-Morrey-0}
\lim_{R \to +\infty} R^{-\frac{\gamma}{p}} \left( \int_{\frac{R}{2} < |x| < R} |f(x)|^p  dx \right)^{\frac{1}{p}} = 0.
\end{equation}
Later, we will describe some important functional spaces that are encompassed by the local Morrey spaces. As mentioned, our main objective in considering these spaces is to study  the Liouville-type problem for the system (\ref{System}) within a fairly general framework.

\medskip

This study involves sharp local estimates on classical solutions to (\ref{System}), often referred to as \emph{Caccioppoli}-type inequalities. When working within the space $M^{p}_{\gamma}(\mathbb{R}^3)$, these estimates hinge on a specific quantity
\begin{equation}\label{q}
    q(p,\gamma):=\frac{\gamma}{p}-\frac{3}{p}+\frac{2}{3},
\end{equation}
which relates the parameters $p$ and $\gamma$. Consequently,  our main result strongly depends  on the sign of this quantity: first of all, we set  $3\leq p < +\infty$ and $0<\gamma<3$ due to some  technical constraints. Thereafter,  in the rectangular region $(p,\gamma)\in [3,+\infty)\times (0,3)$, we examine cases where $q(p,\gamma)<0$, the threshold $q(p,\gamma)=0$, and $q(p,\gamma)>0$. In the  first two cases, when $q(p,\gamma)<0$ and $q(p,\gamma)=0$, the decaying properties given by the spaces $M^{p}_{\gamma}(\Rt)$ and $M^{p}_{\gamma,0}(\Rt)$ are sufficient to prove that classical vanish identically. Conversely,  when $q(p,\gamma)>0$, it is interesting to observe that the  space $M^p_\gamma(\Rt)$ alone may not suffice to resolve the Liouville problem. Additional decay properties are necessary (refer to expression (\ref{Caccioppoli1}) below for further details).

\medskip

On the other hand, it is also interesting to observe the effects of the Rayleigh–Bénard convection term  $\eta \vu \cdot \ve_3$ when studying the Liouville-type problem. Specifically, the aforementioned \emph{Caccioppoli}-type estimates for the system (\ref{System}) rely on beneficial cancellation properties in cross-terms involving the functions $\vu, \vb$, and $\theta$. In this context,  handling with  the Rayleigh–Bénard convection term necessitates the following relationship between the third component of the velocity $\vu$ and the temperature $\theta$:
\begin{equation}\label{Condition1}
u_3  \theta \leq 0, \quad \text{when} \quad \eta > 0,
\end{equation}
for our estimates to hold. Although this requirement is technical, it is not overly restrictive.

\medskip

Returning to our example (\ref{Non-Trivial-Solution}), setting $c=\frac{1}{2}$ reveals that this particular non-trivial solution satisfies (\ref{Condition1}). In conclusion, the key information for solving the Liouville-type problem is consistently provided by the decay properties of solutions.

\medskip

In the case of inactive effects of the Rayleigh–Bénard convection term, when $\eta=0$, our additional assumption (\ref{Condition1}) is not necessary. This scenario is also significant because, as we will elaborate below, when $\eta=0$, the general system (\ref{System}) reduces to other relevant systems in fluid dynamics.

\medskip

Now that we have provided the explanations, we are ready to state our main result. 
\begin{Theoreme}\label{Th-Main} Let $(\vu,\vb,\theta, P)$ be a classical  solution to the coupled system (\ref{System}), and assume that (\ref{Condition1}) holds.  On the other hand, in the local Morrey space $M^{p}_{\gamma}(\Rt)$ defined in (\ref{Local-Morrey}), assume that $(p,\gamma)\in [3,\infty)\times (0,3)$. Moreover, let $q(p,\gamma)$ be the quantity defined in (\ref{q}). Then, the identity 
\[(\vu,\vb,\theta,P)=(0,0,0,0), \]
holds in the following cases: 
\begin{enumerate}
    \item When $\vu, \vb, \theta \in M^{p}_\gamma(\Rt)$ with $q(p,\gamma)<0$.
    \item When $\vu,\vb,\theta  \in M^{p}_{\gamma,0}(\Rt)$  with $q(p,\gamma)=0$.
    \item When $\vu,\vb,\theta \in M^{p}_{\gamma}(\Rt)$,  with $q(p,\gamma)>0$ and with the additional decaying property:
    \begin{equation}\label{Condition2}
        \lim_{R\to +\infty} R^{3q(p,\gamma)-\frac{\gamma}{p}} \left( \int_{\frac{R}{2}<|x|<R} \left( | \vu(x)|^p+ |\vb(x)|^p + | \theta(x)|^p\right) dx \right)^{\frac{1}{p}}=0.
    \end{equation}
\end{enumerate}
\end{Theoreme}

{\bf Remarks and comments:} Comparing points 1 and 2 above, we observe that in the case $q(p,\gamma)<0$, we prove the uniqueness of trivial solutions in the space $M^{p}_{\gamma}(\mathbb{R}^n)$. Meanwhile, in the case $q(p,\gamma)=0$, this result holds in the smaller space $M^{p}_{\gamma,0}(\mathbb{R}^n)$. Nevertheless, both spaces are sufficiently general, and we have the following embeddings in which our Liouville-type result holds. For a brief proof of these facts, see Appendix \ref{AppendixB}. 

\medskip

\emph{The case $q(p,\gamma)<0$}. For the parameters $3<r<\frac{9}{2}$ and  $3\leq p <r$, we have the next chain of embeddings involving the Lebesgue spaces, Lorentz spaces and homogeneous Morrey spaces:  
\begin{equation}\label{Emb-01}
L^r(\Rt) \subset L^{r,\infty}(\Rt)\subset \dot{M}^{r,p}(\Rt)\subset M^p_\gamma(\Rt).    
\end{equation}

\medskip

\emph{The case $q(p,\gamma)=0$.} Setting $p=3$ and $\gamma=1$, for  
a parameter  $\frac{9}{2}<q<+\infty$ we have 
\begin{equation}\label{Emb-02}
 L^{\frac{9}{2}}(\Rt)\subset L^{\frac{9}{2},q}(\Rt)\subset M^{3}_{1,0}(\Rt).   
\end{equation}
These embeddings are of particular interest as they extend the result from \cite[Theorem $X.9.5$]{Galdi}, which was originally established for the Navier-Stokes equations (\ref{NS}), to the coupled setting of the system (\ref{System}). 

\medskip

The case $q(p,\gamma)=0$ also allows us to consider the weighted Lebesgue spaces in which our result holds. For $0<\gamma<3$, we introduce the weight $\ds{w_\gamma(x)=\frac{1}{(1+|x|)^\gamma}}$, and the corresponding weighted Lebesgue space $L^p_{w_\gamma}(\mathbb{R}^3)=L^p(w_\gamma dx)$. Then, for $3\leq p <\frac{9}{2}$, we have
\begin{equation}\label{Emb-03}
L^{p}_{w_\gamma}(\mathbb{R}^3) \subset M^{p}_{\gamma,0}(\mathbb{R}^3).
\end{equation}
Weighted Lebesgue spaces generalize the classical ones in the sense that the weight $w_\gamma(x)$ imposes slower decay properties at infinity. 

\medskip

\emph{The case $q(p,\gamma)>0$}. To our knowledge, due to the strong constraint (\ref{Condition2}), we have not identified any known functional spaces embedded in $M^p_\gamma(\mathbb{R}^n)$ that satisfy (\ref{Condition2}). However, this result highlights the necessity for enhancing the decay properties when $q(p,\gamma)>0$: according to (\ref{Condition2}), the local quantity $\ds{R^{-\frac{\gamma}{p}} \left( \int_{\frac{R}{2}<|x|<R} | f(x)|^p  dx \right)^{\frac{1}{p}}}$ must decay at infinity faster than $R^{3q(p,\gamma)}$. This underscores the eventually sharpness of the condition $q(p,\gamma)\leq 0$.

\medskip

{\bf Applications to related models:} As previously mentioned, the general coupled system (\ref{System}) encompasses other systems of fluid dynamics as special cases. Therefore, Theorem \ref{Th-Main} and the preceding remarks yield new results on the Liouville-type problem for the following systems, which, to the best of our knowledge, have not been studied before. For simplicity, we list them below, although each system possesses significant physical motivation and mathematical interest in its own right.  

\medskip

\begin{itemize}
    \item When $\vb=0$ we have the Bénard system:
\begin{equation*}
    \begin{cases}\vspace{2mm}
    - \Delta \vu + (\vu \cdot \vec{\nabla})\vu +  \vec{\nabla} P -  \theta \ve_3=0, \quad  \text{div}(\vu)=0, \\
    -\Delta \theta + \vu \cdot \vec{\nabla}\theta - \eta \vu \cdot \ve_3=0.   
    \end{cases}
\end{equation*}    
\end{itemize}
The following systems arise from (\ref{System}) with $\eta=0$ (inactive effects of the Rayleigh-Bénard convection term), hence the additional assumption (\ref{Condition1}) is no longer necessary.
\begin{itemize}
    \item When $\eta=0$ we obtain the MHD-Boussinesq system:
    \begin{equation*}
    \begin{cases}\vspace{2mm}
    - \Delta \vu + (\vu \cdot \vec{\nabla})\vu - (\vb \cdot \vec{\nabla}) \vb + \vec{\nabla} \left(P +\frac{1}{2}| \vb|^2 \right) - \theta \ve_3=0, \quad \text{div}(\vu)=0, \\ \vspace{2mm}
 -\Delta \vb + (\vu \cdot \vec{\nabla})\vb - (\vb \cdot \vec{\nabla})\vu =0,  \quad  \text{div}(\vb)=0,\\
  - \Delta \theta + \vu \cdot \vec{\nabla}\theta=0.
    \end{cases}
\end{equation*}
\item When $\eta=0$ and $\vb=0$ we get the Boussinesq system:
\begin{equation*}
    \begin{cases}\vspace{2mm}
    - \Delta \vu + (\vu \cdot \vec{\nabla})\vu +  \vec{\nabla} P -  \theta \ve_3=0, \quad  \text{div}(\vu)=0, \\
    -\Delta \theta + \vu \cdot \vec{\nabla}\theta=0.   
    \end{cases}
\end{equation*}
\item Finally, recall that when $\eta=0$ and $\theta=0$, we recover the MHD equations (\ref{MHD}). In this manner, we establish a result analogous to the one proven in \cite{Li3} for the Hall-MHD equations.
\end{itemize}  

{\bf Strategy of the proof and final comments:} Despite using the same functional framework originally introduced in \cite{Jarrin1}, the coupled system (\ref{System}) presents new technical challenges. As discussed, the active effects of the Rayleigh-Bénard convection term $\eta \vu \cdot \vec{e}_3$ (when $\eta > 0$) necessitate assuming (\ref{Condition1}). It is worth mentioning this assumption can be relaxed by imposing faster decaying properties on the temperature $\theta$:
\begin{equation}\label{Decay-theta}
\lim_{R \to +\infty} R^{2\left( q(p,\gamma) + \frac{5}{6} \right) - \frac{\gamma}{p}} \left( \int_{|x|<R} |\theta(x)|^p  dx \right)^{\frac{1}{p}} = 0,
\end{equation}
even in the case where $q(p,\gamma) \leq 0$. Please refer to Section \ref{Sec:Assumption-Theta}  below for a more  detailed explanation. Nevertheless, this condition on $\theta$ eventually restricts its inclusion in known functional spaces. We have chosen to assume (\ref{Condition1}), which, in light of the example (\ref{Non-Trivial-Solution}), appears to be a more natural \emph{a priori} assumption for classical solutions to (\ref{System}). 

\medskip

On the other hand, the buoyancy force represented by the term $\theta \vec{e}_3$ in the first equation of (\ref{System}) can also pose difficulties when proving the \emph{Caccioppoli}-type estimates. To address this, we circumvent this issue by executing these estimates in a specific sequence. Our general strategy thus focuses initially on the third equation in (\ref{System}) to demonstrate first that $\theta=0$. The remaining identities $(\vu,\vb,P)=(0,0,0)$ then follow from the first and second equations in (\ref{System}).

\section{Proof of Theorem \ref{Th-Main}}
For simplicity, for any $R > 0$, we denote $B_R = \{x \in \mathbb{R}^3 : |x| < R \}$ and $C_R = \{ x \in \mathbb{R}^3 : \frac{R}{2} < |x| < R \}$.
\subsection{\emph{Caccioppoli}-type estimates on $\theta$}
\begin{Proposition}\label{Prop1} Let $(\vu,\vb,\theta,P)$ be a classical solution of (\ref{System}). Moreover,  assume (\ref{Condition1}). Then, for any $R\geq 1$ the following estimate holds:
\begin{equation}\label{Caccioppoli1}
\begin{split}
\int_{B_{\frac{R}{2}}}| \vec{\nabla} \theta  |^2 dx \leq &\, \frac{C}{R^2} \int_{C_R}|\theta|^2 dx + \frac{C}{R}\int_{C_R}|\theta|^2 |\vu| dx,
\end{split}   
\end{equation}
with a constant $C>0$ independent of $R$.
\end{Proposition}
\pv   Our starting point is to introduce a cut-off function as follows: let $\varphi \in \mathcal{C}^{\infty}_0(\Rt)$ be such that $\varphi(x)=1$ for $|x|<\frac{1}{2}$, $\varphi(x)=0$ for $|x|\geq 1$, and $0\leq \varphi(x)\leq 1$ for all $x\in \Rt$. Then, for $R\geq 1$ we define 
\begin{equation}\label{Cut-off}
 \varphi_R(x)=\varphi\left(\frac{x}{R}\right).   
\end{equation}
This function verifies $\text{supp}(\varphi_R)\subset B_R$, for any multi-indice $\alpha \in \mathbb{N}^{3}\setminus \{0\}$ it holds $\text{supp}(\partial^\alpha \varphi_R)\subset C_R$, and we have the estimates:
\begin{equation}\label{Estimates-Cut-off}
 \| \vec{\nabla} \varphi_R\|_{L^\infty}\leq \frac{C}{R}, \quad   \| \Delta \varphi_R \|_{L^\infty}\leq \frac{C}{R^2}, 
\end{equation}
 with a constant $C>0$ independent of $R$.

\medskip

Now, we multiply the third equation in (\ref{System}) by $\varphi_R \theta$ to write
\begin{equation}\label{Iden-theta}
        -\int_{B_R}\Delta \theta \,  \varphi_R \theta dx + \int_{B_R} \vu \cdot \vec{\nabla}\theta \, \varphi_R \theta dx - \eta \int_{B_R}  u_3 \,  \varphi_R \theta dx=0.
\end{equation}
Using integration by parts we explicitly compute each term on the left.  For the first one has
\begin{equation}\label{Estim-Laplacian}
    \begin{split}
  &\,  -\int_{B_R}\Delta \theta   \varphi_R \theta dx= \sum_{i=1}^{3} \int_{B_R} \partial_i \theta \, \partial_i (\varphi_R \theta) dx =   \sum_{i=1}^{3} \int_{C_R} \partial_i \theta (\partial_i \varphi_R) \theta dx  + \sum_{i=1}^{3} \int_{B_R} (\partial_i \theta)^2 \varphi_R dx\\
 =&\,\frac{1}{2}\sum_{i=1}^3 \int_{C_R} \partial_i|\theta|^2 \partial_i \varphi_R dx  +\int_{B_R}\varphi_R |\vec{\nabla}\theta|^2 dx=-\frac{1}{2}\int_{C_R} \Delta \varphi_R |\theta|^2 dx + \int_{B_R}\varphi_R |\vec{\nabla}\theta|^2 dx.   
    \end{split}
\end{equation}
For the second term, since $\text{div}(\vu)=0$, we obtain 
\begin{equation}\label{Estim-non-linear-theta}
  \begin{split}
  &\,  \int_{B_R} \vu \cdot \vec{\nabla}\theta \varphi_R \theta dx=\, \sum_{i=1}^{3}\int_{B_R}u_i (\partial_i \theta) \varphi_R \theta dx = \sum_{i=1}^{3}\int_{B_R} \partial_i (u_i \theta) \varphi_R \theta dx= - \sum_{i=1}^{3} \int_{B_R} u_i \theta \partial_i (\varphi_R \theta) dx\\
  = &\, - \sum_{i=1}^3 \int_{C_R} |\theta|^2 u_i \partial_i \varphi_R dx -  \sum_{i=1}^{3}\int_{B_R} u_i \varphi_R \theta \partial_i \theta dx   = -\int_{C_R}|\theta|^2 (\vu \cdot \vec{\nabla}\varphi_R) dx-\frac{1}{2}\sum_{i=1}^{3}\int_{B_R} u_i  \varphi_R  \partial_i |\theta|^2dx\\
  =&\, -\int_{C_R}|\theta|^2 (\vu \cdot \vec{\nabla}\varphi_R) dx + \frac{1}{2}\sum_{i=1}^{3}\int_{C_R} u_i (\partial_i \varphi_R)|\theta|^2 dx=-\frac{1}{2}  \int_{C_R}|\theta|^2 (\vu \cdot \vec{\nabla}\varphi_R) dx.
  \end{split}  
\end{equation}

With identities (\ref{Estim-Laplacian}) and (\ref{Estim-non-linear-theta}) at hand, we come back to (\ref{Iden-theta}) to get 
\begin{equation}\label{Estim-theta1}
  \int_{B_R}\varphi_R | \vec{\nabla}\theta |^2 dx = \frac{1}{2}\int_{C_R}\Delta \varphi_R |\theta|^2 dx + \frac{1}{2}\int_{C_R}|\theta|^2 (\vu\cdot \vec{\nabla}\varphi_R) dx + \eta \int_{B_R} \theta u_3 \varphi_R dx.
\end{equation}
At this point, we made  the following remarks. On the one hand, recall that $\varphi_R(x)=1$ for $|x|\leq \frac{R}{2}$, hence  
\[ \int_{B_{\frac{R}{2}}} | \vec{\nabla}\theta |^2 dx \leq \int_{B_R}\varphi_R | \vec{\nabla}\theta |^2 dx.\]
On the other hand, recall that  $\varphi_R \geq 0$   and $\eta \geq 0$, consequently,   by our assumption (\ref{Condition1}) (in the case $\eta>0$) we obtain 
\[ \eta \int_{B_R} \theta u_3 \varphi_R dx \leq 0. \]
We thus write
\begin{equation}
  \int_{B_\frac{R}{2}} | \vec{\nabla}\theta |^2 dx \leq  \frac{1}{2}\int_{C_R}\Delta \varphi_R |\theta|^2 dx + \frac{1}{2}\int_{C_R}|\theta|^2 (\vu\cdot \vec{\nabla}\varphi_R) dx,    
\end{equation}
and using the estimates (\ref{Estimates-Cut-off}) we finally obtain the wished estimate (\ref{Caccioppoli1}). Proposition \ref{Prop1}  is proven. \finpv 

\medskip

Estimate (\ref{Caccioppoli1}) will allow us to prove that $\theta=0$. Indeed, since $p\geq 3$ (in particular we have $p>2$) the first term on the right-hand side in (\ref{Caccioppoli1}) can be controlled as follows:
\begin{equation*}
  \frac{C}{R^2} \int_{C_R}|\theta|^2 dx \leq \,  \frac{C}{R^2} R^{6\left(\frac{1}{2}-\frac{1}{p}\right)} \left( \int_{C_R}|\theta|^p dx \right)^{\frac{2}{p}} \leq C\,  R^{6\left(\frac{1}{2}-\frac{1}{p}\right)-2 + \frac{2\gamma}{p}} \left( R^{-\frac{\gamma}{p}} \left(\int_{C_R}|\theta|^p dx \right)^{\frac{1}{p}} \right)^2.
\end{equation*}
Moreover,  recalling the definition of the quantity $q(p,\gamma)$ given in (\ref{q}), we write
\[ 6\left(\frac{1}{2}-\frac{1}{p}\right)-2 + \frac{2\gamma}{p} \leq 2 \left( \frac{\gamma}{p} - \frac{3}{p}+ \frac{1}{2}\right) \leq 2 q(p,\gamma),\]
hence
\begin{equation}\label{Estim-theta-01}
  \frac{C}{R^2} \int_{C_R}|\theta|^2 dx \leq C\, R^{2q(p,\gamma)} \left( R^{-\frac{\gamma}{p}} \left(\int_{C_R}|\theta|^p dx \right)^{\frac{1}{p}} \right)^2.   
\end{equation}

To control the second term on the right-hand side of (\ref{Caccioppoli1}), we use Holder inequalities with the relationship $1=\frac{2}{p}+\frac{1}{r}$. Remark that since $p\geq 3$ we have $r\leq 3 \leq p$. We then write
\begin{equation*}
    \begin{split}
        \frac{C}{R}\int_{C_R}|\theta|^2 |\vu| dx  \leq &\,  \frac{C}{R} \left( \int_{C_R}|\theta|^p dx\right)^{\frac{2}{p}}\left(\int_{C_R}|\vu|^r dx\right)^{\frac{1}{r}}\\
        \leq &\,  \frac{C}{R} \left( \int_{C_R}|\theta|^p dx\right)^{\frac{2}{p}} R^{3\left( \frac{1}{r}-\frac{1}{p}\right)} \left(\int_{C_R}|\vu|^p dx \right)^{\frac{1}{p}}\\
        \leq &\, C\, R^{2-\frac{9}{p}} \left( \int_{C_R}|\theta|^p dx\right)^{\frac{2}{p}} \, \left(\int_{C_R}|\vu|^p dx \right)^{\frac{1}{p}}\\
        \leq &\, C\,  R^{2-\frac{9}{p}-\frac{3\gamma}{p}} \left( R^{-\frac{\gamma}{p}} \left( \int_{C_R}|\theta|^p dx\right)^{\frac{1}{p}}\right)^2 \, \left( R^{-\frac{\gamma}{p}} \left(\int_{C_R}|\vu|^p dx \right)^{\frac{1}{p}} \right)
    \end{split}
\end{equation*}
Here, always by expression (\ref{q})  we have
\[ 2-\frac{9}{p}-\frac{3\gamma}{p}=3\left( \frac{\gamma}{p}-\frac{2}{3}+\frac{3}{p} \right)=3q(p,\gamma), \]
and we obtain
\begin{equation}\label{Estim-theta-02}
      \frac{C}{R}\int_{C_R}|\theta|^2 |\vu| dx  \leq C\, R^{3q(p,\gamma)}  \left( R^{-\frac{\gamma}{p}} \left( \int_{C_R}|\theta|^p dx\right)^{\frac{1}{p}}\right)^2 \, \left( R^{-\frac{\gamma}{p}} \left(\int_{C_R}|\vu|^p dx \right)^{\frac{1}{p}} \right).  
\end{equation}

This way, gathering estimates (\ref{Estim-theta-01}) and (\ref{Estim-theta-02}) in estimate (\ref{Caccioppoli1}), we write our key estimate
\begin{equation}\label{Estim-theta-03}
    \int_{B_{\frac{R}{2}}}|\vec{\nabla}\theta|^2 dx \leq C\left( R^{2q(p,\gamma)}+R^{3q(p,\gamma)} \left( R^{-\frac{\gamma}{p}} \left(\int_{C_R}|\vu|^p dx \right)^{\frac{1}{p}} \right)\right)\left( R^{-\frac{\gamma}{p}} \left( \int_{C_R}|\theta|^p dx\right)^{\frac{1}{p}}\right)^2, 
\end{equation}
hence we must consider the following cases:
\begin{enumerate}
    \item When $q(p,\gamma)<0$. By (\ref{Estim-theta-03}) and the definition of the norm $\| \cdot |_{M^{p}_{\gamma}}$ given in (\ref{Local-Morrey}), we write
    \begin{equation*}
        \int_{B_{\frac{R}{2}}}|\vec{\nabla}\theta|^2 dx \leq C \left( R^{2q(p,\gamma)}+R^{3q(p,\gamma)}\| \vu \|_{M^{p}_{\gamma}} \right)\| \theta\|^{2}_{M^{p}_\gamma}.
    \end{equation*}
Letting $R\to +\infty$ we get that $\vec{\nabla}\theta =0$, and since $\theta \in M^{p}_\gamma(\Rt)$ we conclude that $\theta =0$. 
\item When $q(p,\gamma)=0$.  Coming back to (\ref{Estim-theta-03}) we have 
\begin{equation*}
\int_{B_{\frac{B}{2}}}| \vec{\nabla}\theta|^2 dx \leq C\left(1+\| \vu\|_{M^{p}_{\gamma}}\right)\left(R^{-\frac{\gamma}{p}} \left( \int_{C_R}|\theta|^p dx\right)^{\frac{1}{p}}\right)^2,  
\end{equation*}
and using the information $\theta \in M^{p}_{\gamma,0}(\Rt)$, letting $R\to +\infty$ we get $\theta =0$. 
\item When $q(p,\gamma)>0$. We just rewrite (\ref{Estim-theta-03}) as
\begin{equation*}
\begin{split}
    \int_{B_{\frac{R}{2}}}|\vec{\nabla}\theta|^2 dx \leq &\,  C \left( R^{q(p,\gamma)} \left( R^{-\frac{\gamma}{p}} \left( \int_{C_R}|\theta|^p dx\right)^{\frac{1}{p}}\right)\right)^2\\
    &\, + C\left( R^{q(p,\gamma)} \left( R^{-\frac{\gamma}{p}} \left( \int_{C_R}|\vu|^p dx\right)^{\frac{1}{p}}\right)\right)\left( R^{q(p,\gamma)} \left( R^{-\frac{\gamma}{p}} \left( \int_{C_R}|\theta|^p dx\right)^{\frac{1}{p}}\right)\right)^2,         
\end{split}
\end{equation*}
and using the additional assumption (\ref{Condition2}), letting $R\to +\infty$ it holds $\theta=0$. 
\end{enumerate}

\subsubsection{Removing the assumption (\ref{Condition1})}\label{Sec:Assumption-Theta} As mentioned,  assumption (\ref{Condition1}) can be relaxed by  imposing  faster decay properties on $\theta$ given in (\ref{Decay-theta}). To do this, we  come back to estimate (\ref{Estim-theta1}).  Using H\'older inequalities (with $1=\frac{1}{p}+\frac{1}{p'}$ and $\frac{1}{p'}=\frac{1}{p}+\frac{1}{r}$, hence $1=\frac{2}{p}+\frac{1}{r}$) and the definition of $\varphi_R$ given in (\ref{Cut-off}), in the case $\eta>0$   the term $\ds{\eta \int_{B_R} \theta u_3 \varphi_R \, dx}$ can be estimates as follows:
\begin{equation*}
\begin{split}
 \eta \int_{B_R} \theta u_3 \varphi_R dx \leq &\, \eta \left( \int_{B_R} |\vu|^p dx \right)^{\frac{1}{p}} \left( \int_{B_R}|\varphi_R \, \theta|^{p'} dx \right)^{\frac{1}{p'}} \\
 \leq &\, \eta \left( \int_{B_R} |\vu|^p dx \right)^{\frac{1}{p}}  \left( \int_{B_R} |\theta|^p dx \right)^{\frac{1}{p}}  \left( \int_{B_R} |\varphi_R|^r dx \right)^{\frac{1}{r}} \\
 \leq &\, \eta \left( \int_{B_R} |\vu|^p dx \right)^{\frac{1}{p}}  \left( \int_{B_R} |\theta|^p dx \right)^{\frac{1}{p}} C\, R^{\frac{3}{r}}\\
 = &\, \eta \left( \int_{B_R} |\vu|^p dx \right)^{\frac{1}{p}}  \left( \int_{B_R} |\theta|^p dx \right)^{\frac{1}{p}} C\, R^{3-\frac{6}{p}}.
 \end{split}
\end{equation*}
Then, using the norm $\| \cdot \|_{M^{p}_\gamma}$ defined in (\ref{Local-Morrey}), we get
\begin{equation*}
    \begin{split}
 \eta \int_{B_R} \theta u_3 \varphi_R dx \leq &\, \eta\,  C   \left( \int_{B_R} |\vu|^p dx \right)^{\frac{1}{p}}  R^{3-\frac{6}{p}}\left( \int_{B_R} |\theta|^p dx \right)^{\frac{1}{p}} \\
 \leq &\, \eta\,  C\, R^{-\frac{\gamma}{p}} \left( \int_{B_R} |\vu|^p dx \right)^{\frac{1}{p}} R^{\frac{2\gamma}{p}+3-\frac{6}{p}-\frac{\gamma}{p}} \left( \int_{B_R} |\theta|^p dx \right)^{\frac{1}{p}} \\
 \leq &\, \eta \, C\, \| \vu\|_{M^{p}_{\gamma}} \, R^{\frac{2\gamma}{p}+3-\frac{6}{p}-\frac{\gamma}{p}} \left( \int_{B_R} |\theta|^p dx \right)^{\frac{1}{p}}.
    \end{split}
\end{equation*}
Finally, by the quantity $q(p,\gamma)$ given in (\ref{q}), we write
\[ \frac{2\gamma}{p}+3-\frac{6}{p} = 2\left( \frac{\gamma}{p}-\frac{3}{p}+\frac{3}{2}\right) = 2 \left(q(p,\gamma)+\frac{5}{6}\right). \]
We thus have
\begin{equation*}
\eta \int_{B_R} \theta u_3 \varphi_R dx \leq \eta\,  C\, \| \vu\|_{M^{p}_{\gamma}}\, R^{2 \left(q(p,\gamma)+\frac{5}{6}\right)-\frac{\gamma}{p}}  \left( \int_{B_R} |\theta|^p dx \right)^{\frac{1}{p}}.  
\end{equation*}

This way, for $\eta>0$ estimate (\ref{Caccioppoli1}) becomes:
\begin{equation*}
\begin{split}
\int_{B_{\frac{R}{2}}}| \vec{\nabla} \theta  |^2 dx \leq &\, \frac{C}{R^2} \int_{C_R}|\theta|^2 dx + \frac{C}{R}\int_{C_R}|\theta|^2 |\vu| dx+\eta\, C\, \| \vu\|_{M^{p}_{\gamma}}\, R^{2 \left(q(p,\gamma)+\frac{5}{6}\right)-\frac{\gamma}{p}}  \left( \int_{B_R} |\theta|^p dx \right)^{\frac{1}{p}}.
\end{split}   
\end{equation*}
By this local inequality,  assuming (\ref{Decay-theta}) and following the same computations above, we conclude that $\theta=0$.

\subsection{\emph{Caccioppoli}-type estimates on $\vu$ and $\vb$}

Once we have $\theta =0$, by the third  equation  in  (\ref{System})  we have $u_3=0$ in the case when $\eta>0$. Moreover, the whole system  (\ref{System}) reduces  to the MHD system  (\ref{MHD}). Using this system, we can prove the following Caccioppoli-type estimate involving $\vu$ and $\vb$:
\begin{Proposition}\label{Prop2} Let $(\vu,\vb,P)$ be a smooth solution of (\ref{MHD}). For any $R\geq 1$ the following estimate holds:
\begin{equation}\label{Caccioppoli2}
\int_{B_{\frac{R}{2}}} \left(|\vec{\nabla} \otimes \vu|^2 + |\vec{\nabla}\otimes \vb|^2\right) dx \leq  \frac{C}{R^2} \int_{C_R}\left( |\vu|^2 +|\vb|^2 \right)dx + \frac{C}{R} \int_{C_R}\left(|\vu|^2+|\vb|^2+|P| \right)|\vu|dx.   
\end{equation}
with a constant $C>0$ independent of $R$. 
\end{Proposition}
\pv  We follow very similar ideas in the proof of Proposition \ref{Prop1}. Let consider the cut-off function defined in (\ref{Cut-off}), and 
multiplying  the first equation in (\ref{MHD}) by $\varphi_R \vu$  we get 
\begin{equation}\label{Iden01}
-\int_{B_R}\Delta \vu \cdot \varphi_R \vu dx + \int_{B_R}(\vu\cdot\vec{\nabla})\vu \cdot \varphi_R \vu dx - \int_{B_R}(\vb\cdot\vec{\nabla})\vb \cdot \varphi_R \vu dx + \int_{B_R}\vec{\nabla}\left( P +\frac{1}{2}|\vb|^2\right)\cdot  \varphi_R \vu dx=0.
\end{equation}
The first term on the left is computed as in (\ref{Estim-Laplacian}) and we have
\[ -\int_{B_R}\Delta \vu \cdot \varphi_R \vu dx= -\frac{1}{2}\int_{C_R}\Delta \varphi_R |\vu|^2 dx + \int_{B_R}\varphi_R |\vec{\nabla}\otimes \vu|^2 dx.\]
For the second term, since $\text{div}(\vu)=0$, we write
\begin{equation*}
    \begin{split}
        &\, \int_{B_R} (\vu \cdot \vec{\nabla})\vu \cdot  \varphi_R \vu dx = \sum_{i,j=1}^{3}\int_{B_R} (u_j \partial_j u_i) \varphi_R u_i  dx= \sum_{i,j=1}^{3}\int_{B_R} \partial_j(u_j u_i)\varphi_R u_i  dx\\
        =&\,- \sum_{i,j=1}^{3} \int_{C_R} u^2_i u_j (\partial_j \varphi_R) dx - \sum_{i,j=1}^{3}\int_{B_R} (u_i \partial_j u_i) u_j \varphi_R dx \\
        =&\,  - \int_{C_R} |\vu|^2 (\vu\cdot \vec{\nabla}\varphi_R) dx- \frac{1}{2}\sum_{i,j=1}^{3}\int_{B_R} \partial_j (u^2_i) u_j \varphi_R dx\\
        =&\, - \int_{C_R} |\vu|^2 (\vu\cdot \vec{\nabla}\varphi_R)dx+ \frac{1}{2}\sum_{i,j=1}^{3}\int_{B_R}  u^2_i \partial_j( u_j \varphi_R) dx \\
        =&\, -  \frac{1}{2} \int_{C_R} |\vu|^2 (\vu\cdot \vec{\nabla}\varphi_R)dx. 
    \end{split}
\end{equation*}
The third term follows similar computations, and using the fact that  $\text{div}(\vb)=0$  we obtain
\begin{equation*}
\begin{split}
   &\, - \int_{B_R}(\vb \cdot \vec{\nabla})\vb \cdot \varphi_R \vu dx =- \sum_{i,j=1}^{3}\int_{B_R}\partial_j(b_j b_i)\varphi_R u_i dx \\
   = &\,  \sum_{i,j=1}^{3}\int_{C_R} (u_i b_i) (b_j \partial_j \varphi_R) dx + \sum_{i,j=1}^{3}\int_{B_R} b_j b_i \varphi_R \partial_j u_i dx  \\
   =&\, \int_{C_R} (\vu\cdot \vb)(\vb \cdot \vec{\nabla}\varphi_R) dx + \int_{B_R} (\vb\cdot \vec{\nabla})\vu \cdot \varphi_R \vb dx. 
    \end{split}
\end{equation*}
For the fourth term, using again the divergence-free property of $\vu$, we directly get
\begin{equation*}
    \int_{B_R} \vec{\nabla}\left( P+ \frac{1}{2}|\vb|^2 \right)\cdot \varphi_R \vu dx = - \int_{C_R} \left( P+ \frac{1}{2}|\vb|^2 \right) \vu \cdot \vec{\nabla}\varphi_R dx.
\end{equation*}
This way, identity (\ref{Iden01}) writes down as
\begin{equation}\label{Iden02}
    \begin{split}
 \int_{B_R}\varphi_R |\vec{\nabla}\otimes \vu|^2 dx =&\, \frac{1}{2}\int_{C_R} \Delta \varphi_R |\vu|^2 dx + \frac{1}{2} \int_{C_R}|\vu|^2 (\vu\cdot \vec{\nabla}\varphi_R) dx - \int_{C_R} (\vu\cdot \vb) (\vb \cdot \vec{\nabla}\varphi_R) dx \\
 &\,- \int_{B_R} (\vb\cdot \vec{\nabla})\vu \cdot \varphi_R \vb dx + \int_{C_R}\left( P+\frac{1}{2}|\vb|^2\right) \vu\cdot \vec{\nabla}\varphi_R dx. 
    \end{split}
\end{equation}

Similarly, we multiply the second equation in (\ref{System}) by $\varphi_R \vb$ to get 
\begin{equation*}
\begin{split}
-\int_{B_R} \Delta \vb \cdot \varphi_R \vb dx +\int_{B_R} (\vu\cdot\vec{\nabla})\vb \cdot  \varphi_R \vb dx - \int_{B_R} (\vb\cdot \vec{\nabla})\vu \cdot \varphi_R \vb dx =0.     
\end{split}
\end{equation*}
Here, the first term on the left computes down as
\begin{equation*}
 -\int_{B_R} \Delta \vb \cdot \varphi_R \vb dx= -\frac{1}{2}\int_{C_R}\Delta \varphi_R |\vb|^2 dx + \int_{B_R} \varphi_R |\vec{\nabla}\otimes \vb|^2 dx,   
\end{equation*}
and for the second term we have 
\begin{equation*}
\int_{B_R} (\vu\cdot\vec{\nabla})\vb \cdot  \varphi_R \vb dx = -\frac{1}{2}\int_{C_R}|\vb|^2 (\vu\cdot \vec{\nabla}\varphi_R) dx.
\end{equation*}
Thus, this  identity  rewrites as
\begin{equation}\label{Iden03}
    \begin{split}
 \int_{B_R} \varphi_R |\vec{\nabla}\otimes \vb|^2 dx =&\, \frac{1}{2}\int_{C_R}\Delta \varphi_R |\vb|^2 dx +  \frac{1}{2}\int_{C_R}|\vb|^2 (\vu\cdot \vec{\nabla}\varphi_R) dx+  \int_{B_R} (\vb\cdot \vec{\nabla})\vu \cdot \varphi_R \vb dx.     
    \end{split}
\end{equation}

With identities (\ref{Iden02}) and  (\ref{Iden03})  at hand, and recalling that $\varphi_R(x)=1 $ for $|x|\leq \frac{R}{2}$, we can write 
\begin{equation}\label{Iden05}
  \begin{split}
 & \int_{B_{\frac{R}{2}}} \left( |\vec{\nabla}\otimes \vu|^2 + |\vec{\nabla}\otimes \vb|^2 \right)dx 
\leq    \,  \int_{B_R}\varphi_R \left( |\vec{\nabla}\otimes \vu|^2 + |\vec{\nabla}\otimes \vb|^2  \right)dx\\
    =&\, \frac{1}{2} \int_{C_R}\Delta \varphi_R (|\vu|^2+|\vb|^2)dx + \frac{1}{2}\int_{C_R} (|\vu|^2+|\vb|^2) (\vu \cdot \vec{\nabla}\varphi_R) dx \\
    &\, - \int_{C_R} (\vu\cdot \vb) (\vb \cdot \vec{\nabla}\varphi_R) dx + \int_{C_R}\left( P+\frac{1}{2}|\vb|^2\right) \vu\cdot \vec{\nabla}\varphi_R dx.
  \end{split}  
\end{equation}
Finally, using the estimates (\ref{Estimates-Cut-off}) and rearranging terms we obtain the wished inequality (\ref{Caccioppoli2}). Proposition \ref{Prop2} is proven. \finpv

\subsection{End of the proof of Theorem \ref{Th-Main}}
Once we have estimate (\ref{Caccioppoli2}), the end of the proof of Theorem \ref{Th-Main} follows already known estimates. Indeed, performing the same computations done in (\ref{Estim-theta-01}) one has
\begin{equation}\label{Control-1}
 \frac{C}{R^2} \int_{C_R}\left( |\vu|^2 +|\vb|^2 \right)dx \leq C\, R^{2q(p,\gamma)}\left(R^{-\frac{\gamma}{p}}\left( \int_{C_R} \left( |\vu|^p +|\vb|^p \right) dx \right)^{\frac{1}{p}}\right).   
\end{equation}
Moreover, following the same computations done in (\ref{Estim-theta-02}) (which are essentially based on H\'older inequalities with $1=\frac{2}{p}+\frac{1}{r}$  and $1<r\leq 3 \leq p$) one gets
\begin{equation*}
\begin{split}
&\,\frac{C}{R} \int_{C_R}\left(|\vu|^2+|\vb|^2+|P| \right)|\vu|dx \\
\leq &\,  C\, R^{3q(p,\gamma)} \left( \left( R^{-\frac{\gamma}{p}} \left( \int_{C_R}|\vu|^p dx\right)^{\frac{1}{p}}\right)^2 + \left(R^{-\frac{\gamma}{p}} \left( \int_{C_R}|\vb|^p dx\right)^{\frac{1}{p}}\right)^2 + \left(R^{-\gamma\frac{2}{p}} \left( \int_{C_R}|P|^{\frac{p}{2}} dx\right)^{\frac{2}{p}}\right) \, \right)\times\\
&\, \times \left(R^{-\frac{\gamma}{p}}\left( \int_{C_R} |\vu|^p dx\right)^{\frac{1}{p}}\right).
\end{split}
 \end{equation*}
To control the pressure term $P$, we come back to the first equation in (\ref{MHD}). Applying the divergence operator to each term and since $\text{div}(\vu)=\text{div}(\vb)=0$ we can write
\[ -\Delta \left( P +\frac{1}{2}|\vb|^2 \right)= \text{div}\left( \text{div}(\vu\otimes \vu - \vb \otimes \vb)\right), \]
hence we obtain
\begin{equation}\label{Pressure}
  P= (-\Delta)^{-1}  \text{div}\left( \text{div}(\vu\otimes \vu - \vb \otimes \vb)\right) - \frac{1}{2}|\vb|^2=\sum_{i,j=1}^{3}\mathcal{R}_i \mathcal{R}_j (u_i u_j - b_i b_j) - \frac{1}{2}|\vb|^2,   
\end{equation}
where $\ds{\mathcal{R}_i = \frac{\partial_i}{\sqrt{-\Delta}}}$ is the $i$-st Riesz transform.  At this point, recall that for $0<\gamma<3$ and $1<p<+\infty$  Riesz transforms are bounded in the local Morrey spaces $M^{p}_{\gamma}(\Rt)$ (see [Lemma $2.1$]\cite{Fernandez-Jarrin2}) and we have
\begin{equation*}
\begin{split}
    \left(R^{-\gamma\frac{2}{p}} \left( \int_{C_R}|P|^{\frac{p}{2}} dx\right)^{\frac{2}{p}}\right) \leq &\, \| P\|_{M^{\frac{p}{2}}_{\gamma}} \leq \sum_{i,j=1}^{3} \left\| \mathcal{R}_i \mathcal{R}_j (u_i u_j - b_i b_j) \right\|_{M^{\frac{p}{2}}_{\gamma}}+ \frac{1}{2}\left\| \ |\vb|^2 \right\|_{M^{\frac{p}{2}}_{\gamma}}\\
    \leq &\, C \left( \| \vu\otimes \vu \|_{M^{\frac{p}{2}}_{\gamma}} + \| \vb\otimes \vb \|_{M^{\frac{p}{2}}_{\gamma}}\right) +  \frac{1}{2}\left\| \ |\vb|^2 \right\|_{M^{\frac{p}{2}}_{\gamma}}\\
    \leq &\,  C\left( \| \vu\|^{2}_{M^{p}_{\gamma}}+ \| \vb\|^{2}_{M^{p}_{\gamma}}\right).
    \end{split}
\end{equation*}

Getting back to the inequality above, we find
\begin{equation}\label{Control-2}
\begin{split}
&\,\frac{C}{R} \int_{C_R}\left(|\vu|^2+|\vb|^2+|P| \right)|\vu|dx \\
\leq &\,  C\, R^{3q(p,\gamma)} \left( \left( R^{-\frac{\gamma}{p}} \left( \int_{C_R}|\vu|^p dx\right)^{\frac{1}{p}}\right)^2 + \left(R^{-\frac{\gamma}{p}} \left( \int_{C_R}|\vb|^p dx\right)^{\frac{1}{p}}\right)^2 + C\left( \| \vu\|^{2}_{M^{p}_{\gamma}}+ \| \vb\|^{2}_{M^{p}_{\gamma}}\right) \right)\times\\
&\, \times \left(R^{-\frac{\gamma}{p}}\left( \int_{C_R} |\vu|^p dx \right)^{\frac{1}{p}}\right)\\
\leq &\, C \, R^{3q(p,\gamma)} \left( \| \vu\|^{2}_{M^{p}_{\gamma}}+ \| \vb\|^{2}_{M^{p}_{\gamma}} \right)\left(R^{-\frac{\gamma}{p}}\left( \int_{C_R} |\vu|^p dx \right)^{\frac{1}{p}}\right).
\end{split}
 \end{equation}

Once we dispose of estimates (\ref{Control-1}) and (\ref{Control-2}), we come back to estimate (\ref{Caccioppoli2}) to write
\begin{equation}\label{Estim-Base-MHD}
\begin{split}
 \int_{B_{\frac{R}{2}}} (|\vec{\nabla}\otimes \vu|^2 + |\vec{\nabla}\otimes \vb|^2) dx  \leq &\,   C\, R^{2q(p,\gamma)}\left(R^{-\frac{\gamma}{p}}\left( \int_{C_R}\left( |\vu|^p +|\vb|^p \right) dx \right)^{\frac{1}{p}}\right)\\
 &\, +C\, R^{3q(p,\gamma)} \left( \| \vu\|^{2}_{M^{p}_{\gamma}}+ \| \vb\|^{2}_{M^{p}_{\gamma}} \right)\left(R^{-\frac{\gamma}{p}}\left( \int_{C_R} |\vu|^p dx \right)^{\frac{1}{p}}\right).
\end{split}
\end{equation}
As before, we shall consider the following cases involving the quantity $q(p,\gamma)$ defined in (\ref{q}).
\begin{enumerate}
    \item When $q(p,\gamma)<0$. By (\ref{Estim-Base-MHD}) we have
    \begin{equation*}
 \int_{B_{\frac{R}{2}}} (|\vec{\nabla}\otimes \vu|^2 + |\vec{\nabla}\otimes \vb|^2) dx \leq  C\, R^{2q(p,\gamma)} \left(\| \vu\|_{M^p_\gamma}+\| \vb\|_{M^p_\gamma} \right)  + C\, R^{3q(p,\gamma)}\left( \| \vu\|^2_{M^p_\gamma}+\| \vb\|^2_{M^p_\gamma}\right) \| \vu\|_{M^p_\gamma},
    \end{equation*}
and letting $R\to +\infty$ we get $\vu=0$ and $\vb=0$.
   \item When $q(p,\gamma)=0$.  Always by (\ref{Estim-Base-MHD})  we obtain
  \begin{equation*}
  \begin{split}
     \int_{B_{\frac{R}{2}}} (|\vec{\nabla}\otimes \vu|^2 + |\vec{\nabla}\otimes \vb|^2) dx \leq &\, C  \left(R^{-\frac{\gamma}{p}}\left( \int_{C_R} \left(|\vu|^p +|\vb|^p \right) dx \right)^{\frac{1}{p}}\right) \\
     &\,+ C \left( \| \vu\|^{2}_{M^{p}_{\gamma}}+ \| \vb\|^{2}_{M^{p}_{\gamma}} \right)\left(R^{-\frac{\gamma}{p}}\left( \int_{C_R} |\vu|^p dx \right)^{\frac{1}{p}}\right). 
    \end{split}
  \end{equation*}
In this case we have $\vu, \vb \in M^{p}_{\gamma,0}(\Rt)$ and letting  $R\to +\infty$ we get $\vu=0$ and $\vb=0$.
\item When $q(p,\gamma)>0$. Using again (\ref{Estim-Base-MHD}), since $R\geq 1$  we can write
\begin{equation*}
 \int_{B_{\frac{R}{2}}} (|\vec{\nabla}\otimes \vu|^2 + |\vec{\nabla}\otimes \vb|^2) dx  \leq    C\left(1+\| \vu\|^{2}_{M^{p}_{\gamma}}+ \| \vb\|^{2}_{M^{p}_{\gamma}}\right)\, R^{3q(p,\gamma)}\left(R^{-\frac{\gamma}{p}}\left( \int_{C_R} \left( |\vu|^p +|\vb|^p\right) dx  \right)^{\frac{1}{p}}\right).
 \end{equation*}
By the additional hypothesis (\ref{Condition2}) and letting $R\to +\infty$ we get $\vu=0$ and $\vb=0$.
\end{enumerate}

Finally, once we have the identities $\vu=0$ and $\vb=0$, by expression (\ref{Pressure}) we conclude that $P=0$. Theorem \ref{Th-Main} is proven. \finpv 

\appendix
\section{Appendix}\label{Appendix} Let us consider $(\vu,\vb,\theta,P)$ defined in (\ref{Non-Trivial-Solution}). We will prove that each equation in the system (\ref{System}) is verified by these functions. 

\medskip

\emph{First equation}: we have $\Delta \vu=0$ and since $\vu=\vb$ it directly holds $(\vu \cdot \vec{\nabla})\vu- (\vb \cdot \vec{\nabla})\vb=0$. On the other hand, remark that $\frac{1}{2}|\vb|^2=\frac{a^2}{2}x^{2}_{1}+\frac{b^2}{2}x^{2}_{2}+2c^2 x^{2}_{3}$ and by definition of the pressure $P$ we get $P+\frac{1}{2}|\vb|^2=2c^2 x^{2}_{3}$, hence $\vec{\nabla}\left( P+\frac{1}{2}|\vb|^2=2c^2 x^{2}_{3}\right)=(0,0, 4c^2 x_3)= \theta \ve_3$.  \emph{Second equation}:
this equation directly holds thanks to the identity  $\vu=\vb$ and the fact that $\Delta \vb=0$. \emph{Third equation}: we directly have $\Delta \theta =0$. Moreover, by direct computations we obtain $\vu \cdot \vec{\nabla}\theta=-8c^3 x_3$ and $\vu\cdot \ve_3=-2cx_3$. This way, the identity $\vu \cdot \vec{\nabla}\theta-\vu\cdot \ve_3=0$ is equivalent to the identity $-8c^3 x_3+2cx_3=0$ with holds for any $x_3 \in \R$ as long as $c=\pm \frac{1}{2}$. \emph{Fourth equation}: just observe that  $\vu$ and $\vb$ are divergence-free vector fields as long as $a+b=2c$.

\section{Appendix}\label{AppendixB}
\emph{Proof of embedding (\ref{Emb-01})}. The first three embeddings are well-known, so it is enough to  explain the last one. Just remark that one has $\dot{M}^{r,p}(\Rt)\simeq M^{p}_{\gamma_1}(\Rt)$ with $\gamma_1=3\left(1-\frac{p}{r}\right)$. Here the notation $\simeq$ means that the norms of these spaces are equivalent. Then, for $\gamma_1<\gamma$ we have $M^{p}_{\gamma_1}(\Rt)\subset M^{p}_{\gamma}(\Rt)$. Moreover, the inequality $\gamma_1<\gamma$ is equivalent to $-\frac{3}{r}<\frac{\gamma}{p}-\frac{3}{p}$ and since $q(p,\gamma)<0$ we have the constraint $r<\frac{9}{2}$. 

\medskip

\noindent
\emph{Proof of embedding (\ref{Emb-02})}. Let $f \in L^{\frac{9}{2},q}(\Rt)$. By known embeddings   we have $f\in \dot{M}^{3,\frac{9}{2}}(\Rt)$, and by  the identity $\dot{M}^{3,\frac{9}{2}}(\Rt)= M^{3}_{1}(\Rt)$,  we obtain $f \in M^{3}_{1}(\Rt)$.  Now, let prove that $f \in M^{3}_{1,0}(\Rt)$. For this, we write 
\[  \int_{C_R} \vert f \vert^3 dx= \int_{B_R} \left\vert \mathds{1}_{C_R} f \right\vert^3 dx \leq C \, R \left\Vert  \mathds{1}_{C_R} f  \right\Vert^{3}_{L^{\frac{9}{2},\infty}} \leq C \, R \left\Vert  \mathds{1}_{C_R} f  \right\Vert^{3}_{L^{\frac{9}{2},q}}. \]
By the dominated convergence theorem, which is valid in the space $L^{\frac{9}{2},q}(\Rt)$ for the values $\frac{9}{2}\leq q <+\infty$, see \cite{Chamorro}, we obtain $\ds{\lim_{R \to +\infty} \frac{1}{R} \int_{C_R} \vert f \vert^3\,dx=0}$.

\medskip

\noindent
Finally, embedding (\ref{Emb-03}) was proven in Lemma $2.1$ of our previous work \cite{Fernandez-Jarrin1}.

\end{document}